\documentstyle[12pt,epsfig]{article}
\title{Trilinear relations between the tetrahedron volumes in a cluster}
\author{E.V. Martyushev}
\date{}

\def\be{\begin{eqnarray}}
\def\ee{\end{eqnarray}}
\def\p{\partial}
\def\ov{\overline}

\begin{document}
\maketitle

\begin{abstract}
We prove elegant trilinear formulas connecting products of volumes
of Euclidean tetrahedra with vertices taken from a given set of 6~points.
We propose a way for generalizing those formulas.
\end{abstract}

In this short note we obtain some trilinear relations
between volumes of tetrahedra whose vertices belong to a set
of 6 points, namely $A$, $B$, $C$, $D$, $E$ and $F$ in a 3-dimensional
Euclidean space. This can be imagined as dividing a given octahedron into
tetrahedra. Essentially, these are equalities between determinants
because the volume of a tetrahedron $ABCD$ can be represented as
$$
{1\over 6}\,
\left| \matrix{x_B-x_A&y_B-y_A&z_B-z_A\cr
               x_C-x_A&y_C-y_A&z_C-z_A\cr
               x_D-x_A&y_D-y_A&z_D-z_A}\right|,
$$
where $(x_A,y_A,z_A)$ are the coordinates of vertex $A$ and so on.

It is not claimed that our equalities are new. However, the reader will be
able to see that there are certainly many ``higher analogues'' of our
formulas and maybe {\em some\/} of them will be new. The author will
be grateful for any indications at publications containing formulas of
such kind.

\begin{figure}
\begin{center}
\epsfig{file=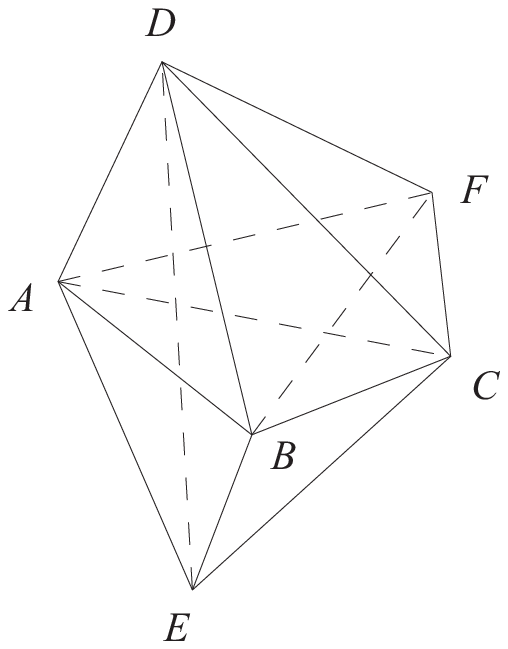}
\end{center}
\caption{}
\label{ris1}
\end{figure}

\begin{figure}
\begin{center}
\epsfig{file=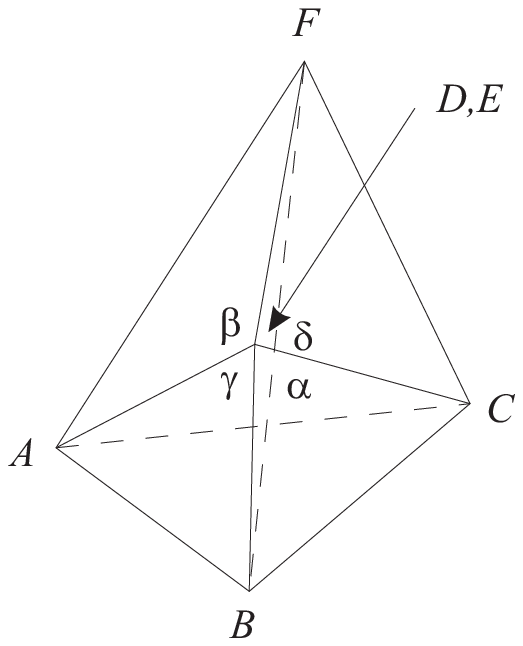}
\end{center}
\caption{}
\label{ris2}
\end{figure}

Consider the Eiclidean octahedron depicted in
Figure~\ref{ris1}.
We will denote the lengths of its edges as
$l_{AB}$, $l_{BC}$ and so on. Draw the diagonal $DE$ and introduce the
following notations for solid angles between planes intersecting at it
as in Figure~\ref{ris2}.
\begin{eqnarray*}
\angle(BDE,CDE)=\alpha,\quad \angle(FDE,ADE)=\beta,\\
\angle(ADE,BDE)=\gamma,\quad \angle(CDE,FDE)=\delta.
\end{eqnarray*}

Then, let us allow the point $D$ to move, but only in such manner that
distances $l_{AD}$ and $l_{FD}$ remain constant. The distances
$l_{BD}$, $l_{CD}$ and $l_{DE}$ will change in some way, and we can
write the following formula:
\be
\label{eq1}
&0 = d(\alpha + \beta + \gamma + \delta) = &\nonumber \\
[1ex]
& = \displaystyle
{\p (\alpha + \beta + \gamma + \delta) \over \p l_{BD}}\,dl_{BD} +
{\p (\alpha + \beta + \gamma + \delta) \over \p l_{CD}}\,dl_{CD} +
{\p (\alpha + \beta + \gamma + \delta) \over \p l_{DE}}\,dl_{DE}.
\ee
The following differential relations can be deduced from simple geometrical
considerations (for details, see~\cite{1}):
\be
\label{eq2}
{\p \alpha \over \p l_{BC}} = {1\over 6}\,{l_{BC}\,l_{DE}
\over V_{\ov{AF}}},&
&{\p \gamma \over \p l_{AB}} = {1\over 6}\,{l_{AB}\,l_{DE}
\over V_{\ov{CF}}},
\ee

\be
\label{eq3}
{\p l_{AB} \over \p l_{BD}} = {l_{BD}\over l_{AB}}\,{V_{\ov{FD}}
\over V_{\ov{AF}}},&
&{\p l_{BC} \over \p l_{CD}} = {l_{CD}\over l_{BC}}\,{V_{\ov{AD}}
\over V_{\ov{AB}}},
\ee

\be
\label{eq4}
{l_{DE}\,dl_{DE} \over V_{\ov{BC}}} = {l_{BD}\,dl_{BD} \over V_{\ov{CE}}} =
{l_{CD}\,dl_{CD} \over V_{\ov{BE}}},
\ee
where, e.g., $V_{\ov{AF}}$ is the volume of tetrahedron $BCDE$, etc.
Then after all substitutions the equality~(\ref{eq1}) reduces to the form
\be
\label{eq5}
{\p (\alpha + \beta + \gamma + \delta) \over \p l_{DE}} =
-{l_{DE}^2 \over 6}\,
\left ({V_{\ov{FD}}\,V_{\ov{CE}} \over V_{\ov{CF}}\,V_{\ov{AF}}\,
V_{\ov{BC}}} +
{V_{\ov{AD}}\,V_{\ov{BE}} \over V_{\ov{AB}}\,V_{\ov{AF}}\,
V_{\ov{BC}}}\right).
\ee

It is important that in obtaining~(\ref{eq5}) we have used the lengths of
$BD$ and $CD$ as ``intermediate'' variables. The use of other pairs of
edges ($CD$ and $FD$, or $FD$ and $AD$, or $AD$ and $BD$)
leads to other relations:
\be
\label{eq6}
{\p (\alpha + \beta + \gamma + \delta) \over \p l_{DE}} =
-{l_{DE}^2 \over 6}\,
\left ({V_{\ov{AD}}\,V_{\ov{FE}} \over V_{\ov{CF}}\,V_{\ov{AF}}\,
V_{\ov{AB}}} +
{V_{\ov{BD}}\,V_{\ov{CE}} \over V_{\ov{AB}}\,V_{\ov{CF}}\,
V_{\ov{BC}}}\right),
\ee
\be
\label{eq7}
{\p (\alpha + \beta + \gamma + \delta) \over \p l_{DE}} =
-{l_{DE}^2 \over 6}\,
\left ({V_{\ov{BD}}\,V_{\ov{AE}} \over V_{\ov{BC}}\,V_{\ov{AF}}\,
V_{\ov{AB}}} +
{V_{\ov{CD}}\,V_{\ov{FE}} \over V_{\ov{AF}}\,V_{\ov{CF}}\,
V_{\ov{BC}}}\right),
\ee
\be
\label{eq8}
{\p (\alpha + \beta + \gamma + \delta) \over \p l_{DE}} =
-{l_{DE}^2 \over 6}\,
\left ({V_{\ov{CD}}\,V_{\ov{BE}} \over V_{\ov{CF}}\,V_{\ov{AB}}\,
V_{\ov{BC}}} +
{V_{\ov{FD}}\,V_{\ov{AE}} \over V_{\ov{AB}}\,V_{\ov{AF}}\,
V_{\ov{CF}}}\right).
\ee
The following three identities obviously follow
 from~(\ref{eq5}),~(\ref{eq6}),~(\ref{eq7}) and~(\ref{eq8}):
\be
\label{eq9}
&V_{\ov{FD}}\,V_{\ov{CE}}\,V_{\ov{AB}} + V_{\ov{AD}}\,V_{\ov{BE}}\,
V_{\ov{CF}}
= 
V_{\ov{AD}}\,V_{\ov{FE}}\,V_{\ov{BC}} + V_{\ov{BD}}\,V_{\ov{CE}}\,
V_{\ov{AF}}
= &\nonumber \\
[1ex]
&V_{\ov{BD}}\,V_{\ov{AE}}\,V_{\ov{CF}} + V_{\ov{CD}}\,V_{\ov{FE}}\,
V_{\ov{AB}}
=
V_{\ov{CD}}\,V_{\ov{BE}}\,V_{\ov{AF}} + V_{\ov{FD}}\,V_{\ov{AE}}\,
V_{\ov{BC}}.
\ee

More equalities of such kind can be obtained by permutations of vertices
and changing the volumes signs, if necessary.

In conclusion, let us note that in considering in a similar way clusters
of tetrahedra with more vertices we must be led to
some ``higher'' multilinear formulae.

\end{document}